\documentclass{amsart}

\usepackage{ adjustbox }
\usepackage{ amsthm }
\usepackage{ amssymb }
\usepackage{ amsmath }
\usepackage{ array }
\usepackage{ booktabs }
\usepackage{ cite }
\usepackage{ color }
\usepackage{ enumitem }
\usepackage{ latexsym }
\usepackage{ url }
\usepackage[all]{ xy }

\theoremstyle{definition}

\newtheorem{definition}{Definition}[section]

\newtheorem{example}[definition]{Example}

\theoremstyle{plain}

\newtheorem{lemma}[definition]{Lemma}

\newtheorem{theorem}[definition]{Theorem}

\makeatletter
\renewcommand{\tocsection}[3]{%
\indentlabel{\@ifnotempty{#2}{\makebox[1.50em][l]{\ignorespaces#1#2.}}}#3}
\renewcommand{\tocsubsection}[3]{%
\indentlabel{\@ifnotempty{#2}{\hspace*{1.50em}\makebox[2.25em][l]{\ignorespaces#1#2.}}}#3}
\renewcommand{\tocsubsubsection}[3]{%
\indentlabel{\@ifnotempty{#2}{\hspace*{3.75em}\makebox[3.00em][l]{\ignorespaces#1#2.}}}#3}
\makeatother

\allowdisplaybreaks

\begin{document}

\title
{Operator identities of multiplicity 3 for associative algebras}

\author{Murray R. Bremner}

\address{Department of Mathematics and Statistics,
University of Saskatchewan,
Saskatoon, Canada}

\email{bremner@math.usask.ca}



\begin{abstract}
We consider algebraic identities for linear operators on associative algebras
in which each term has degree 2 (the number of variables) and multiplicity 3
(the number of occurrences of the operator).
We apply the methods of earlier work by the author and Elgendy
which classified operator identities of degree 2, multiplicities 1 and 2.
We begin with the general operator identity of multiplicity 3 
which has 10 terms and indeterminate coefficients. 
We use the operadic concept of partial composition to generate all consequences
of this identity in degree 3, multiplicity 4.
The coefficient matrix of these consequences has size $105 \times 20$ and
indeterminate entries.
We compute the partial Smith form of this matrix and 
use Gr\"obner bases for determinantal ideals to discover which values of 
the indeterminates produce a matrix of submaximal rank.
The only possible submaximal values of the rank are 16 and 19:
there are 6 new identities of rank 16, and 8 new identities of rank 19.
\end{abstract}

\maketitle

%{\footnotesize\tableofcontents}

%%%%%%%%%%%%%%%%%%%%%%%%%%%%%%%%%%%%%%%%%%%%%%%%%%%%%%%%%%%%%%%%%%%%%%%%%%%%%%%%%%%%%%%%%%%%%%%%
%%%%%%%%%%%%%%%%%%%%%%%%%%%%%%%%%%%%%%%%%%%%%%%%%%%%%%%%%%%%%%%%%%%%%%%%%%%%%%%%%%%%%%%%%%%%%%%%
%%%%%%%%%%%%%%%%%%%%%%%%%%%%%%%%%%%%%%%%%%%%%%%%%%%%%%%%%%%%%%%%%%%%%%%%%%%%%%%%%%%%%%%%%%%%%%%%

\section{Introduction}

%%%%%%%%%%%%%%%%%%%%%%%%%%%%%%%%%%%%%%%%%%%%%%%%%%%%%%%%%%%%%%%%%%%%%%%%%%%%%%%%%%%%%%%%%%%%%%%%

Rota \cite{Rota1995} formulated the problem of classifying algebraic identities satisfied by 
linear operators on associative algebras:
``I have posed the problem of finding all possible algebraic identities
that can be satisfied by a linear operator on an [associative] algebra.
Simple computations show that the possibilities are very few,
and the problem of classifying all such identities is very probably completely solvable.''

For a list of well-known homogeneous and inhomogeneous operator identities, 
see \cite[Table 1]{BE2022}.
That paper also introduces some previously unknown operator identities
\cite[Theorem 6.12]{BE2022}.
For further information about the classification of operator identities,
see \cite{GaoGuo2017,Guo2009,GSZ1,GSZ2}
and the other papers cited in \cite[References]{BE2022}.

A more recent paper that deserves mention is Wang et al.~\cite{WZG2022} which proves
that all the operator identities in \cite[Theorem 6.12]{BE2022} satisfy 
the Gr\"obner-Shirshov property.
A new direction in the study of Rota's problem is the classification of 
algebraic identities for linear operators on Lie algebras,
see Zhang et al.~\cite{ZGG2023}.
In particular, the Lie analogues of the identities in \cite[Theorem 6.12]{BE2022}
(obtained by replacing each associative product by the Lie bracket) have been shown
to satisfy the Gr\"obner-Shirshov property by Zhang et al.~\cite{ZGWF2024}.
It seems to be an open problem to classify operator identities for other varieties
of algebras.

%%%%%%%%%%%%%%%%%%%%%%%%%%%%%%%%%%%%%%%%%%%%%%%%%%%%%%%%%%%%%%%%%%%%%%%%%%%%%%%%%%%%%%%%%%%%%%%%
%%%%%%%%%%%%%%%%%%%%%%%%%%%%%%%%%%%%%%%%%%%%%%%%%%%%%%%%%%%%%%%%%%%%%%%%%%%%%%%%%%%%%%%%%%%%%%%%
%%%%%%%%%%%%%%%%%%%%%%%%%%%%%%%%%%%%%%%%%%%%%%%%%%%%%%%%%%%%%%%%%%%%%%%%%%%%%%%%%%%%%%%%%%%%%%%%

\section{Preliminaries}

%%%%%%%%%%%%%%%%%%%%%%%%%%%%%%%%%%%%%%%%%%%%%%%%%%%%%%%%%%%%%%%%%%%%%%%%%%%%%%%%%%%%%%%%%%%%%%%%

\subsection{Basic definitions}

An \emph{associative algebra} is a vector space $A$ with a bilinear operation
$B\colon A \times A \to A$ which is denoted by juxtaposition $(a,b) \mapsto B(a,b) = ab$
and which satisfies the associativity relation $(ab)c = a(bc)$ for all $a, b, c \in A$.
A \emph{monomial of degree $p$} in $A$ is a product $a_1 a_2 \cdots a_p$ where the $p$ indeterminates
represent arbitrary elements of $A$.
We write $L\colon A \to A$ for a linear operator on the underlying vector space of $A$.
An \emph{operator monomial of degree $p$ and multiplicity $q$} is a monomial of degree $p$
into which are inserted $q$ occurrences of $L$.
For example, the 10 monomials of degree 2 and multiplicity 3 are as follows,
where we use the asterisk as a generic argument symbol,
since the underlying structure is a nonsymmetric operad in which each monomial
contains the identity permutation of the variables:
\begin{equation}
\label{monomials23}
\begin{array}{lllll}
L^3({\ast}{\ast}) 
&\quad  L^2( L({\ast}) {\ast} ) 
&\quad  L( L^2({\ast}) {\ast} ) 
&\quad  L^3({\ast}) {\ast} 
&\quad  L^2( {\ast} L({\ast}) ) 
\\[2pt]
L( L({\ast}) L({\ast}) )  
&\quad  L^2({\ast}) L({\ast}) 
&\quad  L({\ast} L^2({\ast})) 
&\quad  L({\ast}) L^2({\ast}) 
&\quad  {\ast} L^3({\ast})
\end{array}
\end{equation}
A (homogeneous) \emph{operator polynomial} $I$ is a linear combination of operator monomials 
of the same degree and multiplicity.
An \emph{operator identity} is an equation of the form $I \equiv 0$ which holds 
for all values of the indeterminates.

%%%%%%%%%%%%%%%%%%%%%%%%%%%%%%%%%%%%%%%%%%%%%%%%%%%%%%%%%%%%%%%%%%%%%%%%%%%%%%%%%%%%%%%%%%%%%%%%

\subsection{Methods and results}

The author and Elgendy \cite{BE2022} introduced a new approach to 
the classification of operator identities,
using computational linear algebra over polynomial rings 
and computational commutative algebra (especially Gr\"obner bases for determinantal ideals).
We quote the outline of this method from that paper:
``We start with an operator identity with indeterminate coefficients.
We apply the operadic concept of partial composition
to produce a set of new identities for which both the degree and multiplicity
increase by 1.
We call these new identities the \emph{consequences} of the original identity.
We store the coefficient vectors of these consequences [as the column vectors]
in a matrix over the polynomial ring generated by the indeterminates.
We call this the \emph{matrix of consequences} of the original identity.
We introduce what we call the \emph{rank principle}, namely that significant operator identities
correspond to those values of the coefficients which produce submaximal rank
of the matrix of consequences.''
We use Gr\"obner bases for determinantal ideals to determine how the rank
of the matrix of consequences depends on the values of the indeterminates.
We find that extremely few values of the coefficients produce a matrix of submaximal rank.

All of these computations were performed with the computer algebra system Maple
\cite{Maple}.
We assume the base field $\mathbb{F}$ has characteristic 0.

We consider only operator identities which are homogeneous of degree 2
and multiplicity 3.
Our main result is Theorem \ref{theorem23}:
we obtain 6 identities for rank 16, and 8 identities for rank 19.
For rank 16, we obtain the derivation identity for $L^3$ together with
some new identities involving one or two terms.
For rank 19, we obtain identities generalizing the left and right versions
of identity $A$ from \cite[Theorem 6.12]{BE2022}, 
some of which have the imaginary unit $i$ among their coefficients.

%%%%%%%%%%%%%%%%%%%%%%%%%%%%%%%%%%%%%%%%%%%%%%%%%%%%%%%%%%%%%%%%%%%%%%%%%%%%%%%%%%%%%%%%%%%%%%%%

\subsection{Narayana numbers and Dyck words}

The \emph{Narayana numbers} \cite[A001263]{OEIS} are defined by
\[
N(i,j) = \frac1i \binom{i}{j} \binom{i}{j{-}1}
\qquad
(i \ge j \ge 1).
\]
We quote from \cite{BE2022}:
``For $i \ge 1$, consider strings of length $2i$ consisting of
$i$ left parentheses and $i$ right parentheses.
Such a string is \emph{balanced} if in every initial substring
the number of right parentheses is no greater than the number of left parentheses.
A \emph{nesting} is a substring of the form $()$.''
Such a string of balanced parentheses is called a \emph{Dyck word}.

\begin{lemma}
\label{narayanalemma}
$N(i,j)$ is the number of Dyck words of length $2i$ which contain $j$ nestings.
\end{lemma}

Let $\mathcal{M}(p,q)$ denote the set of all operator monomials of
degree $p$ and multiplicity $q$,
and let $\mathcal{O}(p,q)$ be the vector space with basis $\mathcal{M}(p,q)$.
Let $\mathcal{S}(p,q)$ be the set of all Dyck words
containing $p{+}q$ balanced pairs of parentheses, $p$ of which are nestings.
We define a bijection from $\mathcal{S}(p,q)$ to $\mathcal{M}(p,q)$.
Given an element of $\mathcal{S}(p,q)$, we replace each nesting by the argument symbol $\ast$.
We then insert the operator symbol $L$
immediately to the left of each remaining left parenthesis,
and leave the remaining right parentheses unchanged.
This bijection induces
a \emph{lexicographic order} on operator monomials $v, w \in \mathcal{M}(p,q)$.
Replace $v$ and $w$ by the corresponding Dyck words $s(v)$ and $s(w)$.
Find the leftmost position in which $s(v)$ and $s(w)$ differ, 
and assume that $($ precedes $)$.
We used this order in display (\ref{monomials23}) and Table \ref{monomials34}.
Hence for $p \ge 1$ and $q \ge 0$ we have
\[
\dim \mathcal{O}(p,q) = N(p{+}q,p) = \frac{1}{p{+}q} \binom{p{+}q}{p} \binom{p{+}q}{p{-}1}.
\]

%%%%%%%%%%%%%%%%%%%%%%%%%%%%%%%%%%%%%%%%%%%%%%%%%%%%%%%%%%%%%%%%%%%%%%%%%%%%%%%%%%%%%%%%%%%%%%%%

\begin{table}
$
\begin{array}{llll}
L^4({*}{*}{*})  &\,  %  ((((()()()))))    1 
L^3(L({*}{*}){*})  &\,  %  ((((()())())))    2 
L^2(L^2({*}{*}){*})  &\,  %  ((((()()))()))    3 
L(L^3({*}{*}){*})   %  ((((()())))())    4 
\\
L^4({*}{*}){*}  &\,   %  ((((()()))))()    5 
L^3(L({*}){*}{*})  &\,  %  ((((())()())))    6 
L^2(L(L({*}){*}){*})  &\,  %  ((((())())()))    7 
L(L^2(L({*}){*}){*})    %  ((((())()))())    8 
\\
L^3(L({*}){*}){*}  &\,  %  ((((())())))()    9 
L^2(L^2({*}){*}{*})  &\,   %  ((((()))()()))   10 
L(L(L^2({*}){*}){*})  &\,  %  ((((()))())())   11 
L^2(L^2({*}){*}){*}    %  ((((()))()))()   12 
\\
L(L^3({*}){*}{*})  &\,  %  ((((())))()())   13 
L(L^3({*}){*}){*}  &\,  %  ((((())))())()   14 
L^4({*}){*}{*}  &\,  %  ((((()))))()()   15 
L^3({*}L({*}{*}))    %  (((()(()()))))   16 
\\
L^3({*}L({*}){*})  &\,  %  (((()(())())))   17 
L^2(L({*}L({*})){*})  &\,  %  (((()(()))()))   18 
L(L^2({*}L({*})){*})  &\,  %  (((()(())))())   19 
L^3({*}L({*})){*}     %  (((()(()))))()   20 
\\
L^3({*}{*}L({*}))  &\,  %  (((()()(()))))   21 
L^2(L({*}{*})L({*}))  &\,  %  (((()())(())))   22 
L(L^2({*}{*})L({*}))  &\,  %  (((()()))(()))   23 
L^3({*}{*})L({*})    %  (((()())))(())   24 
\\
L^2(L({*})L({*}{*}))   &\,  %  (((())(()())))   25 
L^2(L({*})L({*}){*})  &\,  %  (((())(())()))   26 
L(L(L({*})L({*})){*})  &\,  %  (((())(()))())   27 
L^2(L({*})L({*})){*}    %  (((())(())))()   28 
\\
L^2(L({*}){*}L({*}))  &\,  %  (((())()(())))   29 
L(L(L({*}){*})L({*})) &\,   %  (((())())(()))   30 
L^2(L({*}){*})L({*})  &\,  %  (((())()))(())   31 
L(L^2({*})L({*}{*}))     %  (((()))(()()))   32 
\\
L(L^2({*})L({*}){*})  &\,  %  (((()))(())())   33 
L(L^2({*})L({*})){*}  &\,  %  (((()))(()))()   34 
L(L^2({*}){*}L({*}))  &\,  %  (((()))()(()))   35 
L(L^2({*}){*})L({*})    %  (((()))())(())   36 
\\
L^3({*})L({*}{*})  &\,  %  (((())))(()())   37 
L^3({*})L({*}){*}  &\,  %  (((())))(())()   38 
L^3({*}){*}L({*})  &\,  %  (((())))()(())   39 
L^2({*}L^2({*}{*}))     %  ((()((()()))))   40 
\\
L^2({*}L(L({*}){*}))  &\,  %  ((()((())())))   41 
L^2({*}L^2({*}){*})  &\,  %  ((()((()))()))   42 
L(L({*}L^2({*})){*})  &\,  %  ((()((())))())   43 
L^2({*}L^2({*})){*}    %  ((()((()))))()   44 
\\
L^2({*}L({*}L({*})))   &\,  %  ((()(()(()))))   45 
L^2({*}L({*})L({*}))  &\,  %  ((()(())(())))   46 
L(L({*}L({*}))L({*}))  &\,  %  ((()(()))(()))   47 
L^2({*}L({*}))L({*})    %  ((()(())))(())   48 
\\
L^2({*}{*}L^2({*}))  &\,  %  ((()()((()))))   49 
L(L({*}{*})L^2({*}))   &\,  %  ((()())((())))   50 
L^2({*}{*})L^2({*})  &\,  %  ((()()))((()))   51 
L(L({*})L^2({*}{*}))    %  ((())((()())))   52 
\\
L(L({*})L(L({*}){*}))  &\,  %  ((())((())()))   53 
L(L({*})L^2({*}){*})  &\,  %  ((())((()))())   54 
L(L({*})L^2({*})){*}   &\,  %  ((())((())))()   55 
L(L({*})L({*}L({*})))    %  ((())(()(())))   56 
\\
L(L({*})L({*})L({*}))  &\,  %  ((())(())(()))   57 
L(L({*})L({*}))L({*})  &\,  %  ((())(()))(())   58 
L(L({*}){*}L^2({*}))  &\,  %  ((())()((())))   59 
L(L({*}){*})L^2({*})     %  ((())())((()))   60 
\\
L^2({*})L^2({*}{*})  &\,  %  ((()))((()()))   61 
L^2({*})L(L({*}){*})  &\,  %  ((()))((())())   62 
L^2({*})L^2({*}){*}  &\,  %  ((()))((()))()   63 
L^2({*})L({*}L({*}))    %  ((()))(()(()))   64 
\\
L^2({*})L({*})L({*})  &\,   %  ((()))(())(())   65 
L^2({*}){*}L^2({*})  &\,  %  ((()))()((()))   66 
L({*}L^3({*}{*}))  &\,  %  (()(((()()))))   67 
L({*}L^2(L({*}){*}))    %  (()(((())())))   68 
\\
L({*}L(L^2({*}){*}))  &\,  %  (()(((()))()))   69 
L({*}L^3({*}){*})   &\,  %  (()(((())))())   70 
L({*}L^3({*})){*}  &\,  %  (()(((()))))()   71 
L({*}L^2({*}L({*})))    %  (()((()(()))))   72 
\\
L({*}L(L({*})L({*})))  &\,  %  (()((())(())))   73 
L({*}L^2({*})L({*}))  &\,  %  (()((()))(()))   74 
L({*}L^2({*}))L({*})   &\,  %  (()((())))(())   75 
L({*}L({*}L^2({*})))    %  (()(()((()))))   76 
\\
L({*}L({*})L^2({*}))  &\,  %  (()(())((())))   77 
L({*}L({*}))L^2({*})  &\,  %  (()(()))((()))   78 
L({*}{*}L^3({*}))  &\,  %  (()()(((()))))   79 
L({*}{*})L^3({*})     %  (()())(((())))   80 
\\
L({*})L^3({*}{*})  &\,  %  (())(((()())))   81 
L({*})L^2(L({*}){*})  &\,  %  (())(((())()))   82 
L({*})L(L^2({*}){*})  &\,  %  (())(((()))())   83 
L({*})L^3({*}){*}    %  (())(((())))()   84 
\\
L({*})L^2({*}L({*}))  &\,   %  (())((()(())))   85 
L({*})L(L({*})L({*}))  &\,  %  (())((())(()))   86 
L({*})L^2({*})L({*})  &\,  %  (())((()))(())   87 
L({*})L({*}L^2({*}))    %  (())(()((())))   88 
\\
L({*})L({*})L^2({*})  &\,  %  (())(())((()))   89 
L({*}){*}L^3({*})   &\,  %  (())()(((())))   90 
{*}L^4({*}{*})  &\,  %  ()((((()()))))   91 
{*}L(L^3({*}){*})    %  ()((((())())))   92 
\\
{*}L^2(L^2({*}){*})  &\,  %  ()((((()))()))   93 
{*}L(L^3({*}){*})  &\,  %  ()((((())))())   94 
{*}L^4({*}){*}   &\,  %  ()((((()))))()   95 
{*}L^3({*}L({*}))    %  ()(((()(()))))   96 
\\
{*}L^2(L({*})L({*}))  &\,  %  ()(((())(())))   97 
{*}L(L^2({*})L({*}))  &\,  %  ()(((()))(()))   98 
{*}L^3({*})L({*})  &\,  %  ()(((())))(())   99 
{*}L^2({*}L^2({*}))     %  ()((()((()))))  100 
\\
{*}L(L({*})L^2({*}))  &\,  %  ()((())((())))  101 
{*}L^2({*})L^2({*})  &\,  %  ()((()))((()))  102 
{*}L({*}L^3({*}))  &\,  %  ()(()(((()))))  103 
{*}L({*})L^3({*})    %  ()(())(((())))  104 
\\
{*}{*}L^4({*})     %  ()()((((()))))  105 
\end{array}
$
\medskip
\caption{
The 105 operator monomials of degree 3, multiplicity 4
which form the lexicographically ordered basis $\mathcal{M}(3,4)$ of $\mathcal{O}(3,4)$}
\label{monomials34}
\end{table}

%%%%%%%%%%%%%%%%%%%%%%%%%%%%%%%%%%%%%%%%%%%%%%%%%%%%%%%%%%%%%%%%%%%%%%%%%%%%%%%%%%%%%%%%%%%%%%%%

\subsection{Partial compositions}

We refer to the author and Dotsenko \cite{BD-book} for basic information on algebraic operads.
Partial compositions are defined as follows for $m \in \mathcal{M}(p,q)$,
where $B$ denotes the binary associative operation and $L$ denotes the unary operator:
\begin{enumerate}
\item
For $1 \le i \le p$, we define $m \circ_i B$ as follows.
Replace the $i$-th argument symbol $\ast$ by the product $\ast\ast$,
obtaining $m \circ_i B \in \mathcal{M}(p{+}1,q)$.
\item
For $1 \le j \le 2$, we define $B \circ_j m$ as follows.
Multiply $m$ on the right ($j = 1$) or the left ($j = 2$) by
the argument symbol $\ast$, obtaining
$B \circ_j m \in \mathcal{M}(p{+}1,q)$.
\item
For $1 \le i \le p$, we define $m \circ_i L$ as follows:
apply the linear operator $L$ to the $i$-th argument symbol,
obtaining $m \circ_i L \in \mathcal{M}(p,q{+}1)$.
\item
To define $L \circ m$, apply $L$ to $m$< obtaining $L \circ m \in \mathcal{M}(p,q{+}1)$.
\end{enumerate}
Partial composition extends linearly from $\mathcal{M}(p,q)$ to $\mathcal{O}(p,q)$.

Let $R \in \mathcal{O}(p,q)$.
By a \emph{consequence} of $R$ we mean an element of $\mathcal{O}(p+1,q+1)$
obtained by applying partial compositions to $R$:
we either increase the degree first and then the multiplicity,
or the multiplicity first and then the degree.

\begin{example}
We use the operator monomials in display \eqref{monomials23}.
For each of these monomials $m$ we will calculate the partial composition
$(m \circ_1 B) \circ_2 L$, 
increasing the degree first and then the multiplicity.
The resulting monomials appear in Table \ref{monomials34},
which displays all operator monomials of degree 3 and multiplicity 4.
Since partial composition is a linear map, this shows how to compute
the consequence $( R \circ_1 B ) \circ_2 L$ of the general linear combination $R$ 
of the monomials \eqref{monomials23}:
\begin{align*}
( L^3({\ast}{\ast}) \circ_1 B ) \circ_2 L
&=
L^3({\ast}{\ast}{\ast}) \circ_2 L
=
L^3({\ast}L({\ast}){\ast}),
\\
( L^2( L({\ast}) {\ast} ) \circ_1 B ) \circ_2 L 
&=
L^2( L({\ast}{\ast}) {\ast} ) \circ_2 L 
=
L^2( L({\ast}L({\ast})) {\ast} ),
\\
( L( L^2({\ast}) {\ast} ) \circ_1 B ) \circ_2 L 
&=
L( L^2({\ast}{\ast}) {\ast} ) \circ_2 L 
=
L( L^2({\ast}L({\ast})) {\ast} ),
\\
( L^3({\ast}) {\ast} \circ_1 B ) \circ_2 L 
&=
L^3({\ast}{\ast}) {\ast} \circ_2 L 
=
L^3({\ast}L({\ast})) {\ast},
\\
( L^2( {\ast} L({\ast}) ) \circ_1 B ) \circ_2 L 
&=
L^2( {\ast}{\ast} L({\ast}) ) \circ_2 L 
=
L^2( {\ast}L({\ast}) L({\ast}) ),
\\
( L( L({\ast}) L({\ast}) ) \circ_1 B ) \circ_2 L  
&=
L( L({\ast}{\ast}) L({\ast}) ) \circ_2 L  
=
L( L({\ast}L({\ast})) L({\ast}) ),
\\
( L^2({\ast}) L({\ast}) \circ_1 B ) \circ_2 L
&=
L^2({\ast}{\ast}) L({\ast}) \circ_2 L
=
L^2({\ast}L({\ast})) L({\ast}),
\\
( L({\ast} L^2({\ast})) \circ_1 B ) \circ_2 L 
&=
L({\ast}{\ast} L^2({\ast})) \circ_2 L 
=
L( {\ast} L({\ast}) L^2({\ast}) ),
\\
( L({\ast}) L^2({\ast}) \circ_1 B ) \circ_2 L 
&=
L({\ast}{\ast}) L^2({\ast}) \circ_2 L 
=
L({\ast}L({\ast})) L^2({\ast}),
\\
( {\ast} L^3({\ast}) \circ_1 B ) \circ_2 L
&=
{\ast}{\ast} L^3({\ast}) \circ_2 L
=
{\ast} L({\ast}) L^3({\ast}).
\end{align*}
\end{example}

%%%%%%%%%%%%%%%%%%%%%%%%%%%%%%%%%%%%%%%%%%%%%%%%%%%%%%%%%%%%%%%%%%%%%%%%%%%%%%%%%%%%%%%%%%%%%%%%

\subsection{Linear algebra and commutative algebra}

Reference \cite[Chapters 7-10]{BD-book} discusses applications of
computational commutative algebra to classification of operads,
and similar methods will be applied in this paper.

Let $R \in \mathcal{O}(p,q)$.
We write $M(R)$ for the \emph{matrix of consequences} of $R$,
which has the coefficient vectors of the consequences as its column vectors.
Thus $M(R)$ has $N(p{+}q{+}2,p{+}1)$ rows, and its number of columns equals
the number of (distinct) consequences of $R$.
The $(i,j)$ entry of $M(R)$ is the
coefficient of the $i$-th monomial in the ordered basis $\mathcal{M}(p{+}1,q{+}1)$
in the $j$-th consequence of $R$.
If $R$ has indeterminate coefficients then
$M(R)$ is a matrix over the polynomial ring in $N(p{+}q,p)$ variables.

We reformulate Rota's problem in terms of linear algebra over polynomial rings.
Let $R \in \mathcal{O}(p,q)$ with indeterminate coefficients.
Our goal is to understand how the rank of $M(R)$ depends
on the values of the indeterminates.
We find that very few values of rank$(M(R))$ are possible,
and that each possible submaximal rank corresponds to a small set of coefficient vectors.

We rely on the following well-known facts.
Let $M$ be an $m \times n$ matrix with coefficients in $\mathbb{F}$.
For $0 \le r \le \min(m,n)$ we have $\mathrm{rank}(M) = r$ if and only if
\begin{enumerate}
\item[(i)]
some $r \times r$ submatrix of $M$ has nonzero determinant, and
\item[(ii)]
every $(r{+}1) \times (r{+}1)$ submatrix of $M$ has determinant 0.
\end{enumerate}
Let $M$ be an $m \times n$ matrix with coefficients in $\mathbb{F}[x_1,\dots,x_n]$.
Let $I(M,r)$ be the $r$-th \emph{determinantal ideal} of $M$ which is
generated by the determinants of all $r \times r$ submatrices of $M$.
Let $Z(M,r) \subseteq \mathbb{F}^k$ be the zero set of $I(M,r)$.
For $0 \le r < s \le \min(m,n)$ we have $I(M,r) \supseteq I(M,s)$
and hence $Z(M,r) \subseteq Z(M,s)$.
For $X = (a_1,\dots,a_k) \in \mathbb{F}^k$ we obtain $M_X$ by substituting 
$x_i = a_i$ ($i=1,\dots,k$) in $M$.
Then $M_X$ has rank $r$ if and only if $X \in Z(M,r{+}1)$ but $X \notin Z(M,r)$.

%%%%%%%%%%%%%%%%%%%%%%%%%%%%%%%%%%%%%%%%%%%%%%%%%%%%%%%%%%%%%%%%%%%%%%%%%%%%%%%%%%%%%%%%%%%%%%%%
%%%%%%%%%%%%%%%%%%%%%%%%%%%%%%%%%%%%%%%%%%%%%%%%%%%%%%%%%%%%%%%%%%%%%%%%%%%%%%%%%%%%%%%%%%%%%%%%
%%%%%%%%%%%%%%%%%%%%%%%%%%%%%%%%%%%%%%%%%%%%%%%%%%%%%%%%%%%%%%%%%%%%%%%%%%%%%%%%%%%%%%%%%%%%%%%%

\section{Operator identities of degree 2 and multiplicity 3}
\label{deg2mul3}

%%%%%%%%%%%%%%%%%%%%%%%%%%%%%%%%%%%%%%%%%%%%%%%%%%%%%%%%%%%%%%%%%%%%%%%%%%%%%%%%%%%%%%%%%%%%%%%%

\subsection{The matrix of consequences}
\label{construct23}

The general operator identity $R \in \mathcal{O}(2,3)$ is
\begin{align*}
&
a_1 L^3({\ast}{\ast}) +
a_2 L^2( L({\ast}) {\ast} ) +
a_3 L( L^2({\ast}) {\ast} ) +
a_4 L^3({\ast}) {\ast} +
a_5 L^2( {\ast} L({\ast}) ) + {}
\\
&
a_6 L( L({\ast}) L({\ast}) ) + 
a_7 L^2({\ast}) L({\ast}) +
a_8 L({\ast} L^2({\ast})) +
a_9 L({\ast}) L^2({\ast}) +
a_{10} {\ast} L^3({\ast}),
\end{align*}
where $a_1, \dots, a_{10}$ are indeterminate parameters from $\mathbb{F}$.
Without loss of generality, we may scale $R$ so that its first nonzero coefficient is 1.
Thus we may split the problem into 10 mutually exclusive cases:
for $1 \le k \le 10$, case $k$ considers the matrix of consequences 
in which coefficients $a_1,\dots,a_{k-1}$ are 0,
$a_k = 1$, and $a_{k+1},\dots,a_{10}$ are free parameters.

We have altogether 20 distinct consequences of $R$ in degree 3, multiplicity 4.
The monomials in the consequences belong to the ordered monomial basis
of $\mathcal{O}(3,4)$ given in Table \ref{monomials34}.
The matrix of consequences $M(R)$ has size $105 \times 20$.
Our goal is to understand how the rank of this matrix depends on the values of
the parameters.

We note that rows 42, 57, 66 of the matrix of consequences are 0 corresponding to monomials
$L^2({*}L^2({*}){*})$,
$L(L({*})L({*})L({*}))$,
$L^2({*}){*}L^2({*})$.
The reason is that these are the monomials of degree 3, multiplicity 4 
which cannot be obtained by partial compositions
from monomials of degree 2, multiplicity 3.

Throughout the computations described below we use the deglex monomial order
with $a_1 \prec a_2 \prec \cdots \prec a_{10}$.

%%%%%%%%%%%%%%%%%%%%%%%%%%%%%%%%%%%%%%%%%%%%%%%%%%%%%%%%%%%%%%%%%%%%%%%%%%%%%%%%%%%%%%%%%%%%%%%%

\subsection{Case 1: $a_1 = 1$, and $a_2, \dots, a_{10}$ are free}

We substitute $a = 1$ into the matrix $M$ of consequences.
We compute the partial Smith form of the resulting matrix and obtain
\[
M
\xrightarrow{\quad\text{PSF}\quad}
\left[
\begin{array}{cc}
I_{16} &\!\! O_{16,4} \\
O_{89,16} &\!\! C
\end{array}
\right]
\]
where $I$ and $O$ denote the identity and zero matrices of the given sizes, and
the lower right block $C$ has size $89 \times 4$.
We delete the four zero rows in $C$ leaving a $85 \times 4$ matrix also called $C$.

The partial Smith form is obtained by using elementary row operations 
to repeatedly swap matrix entries equal to 1 onto the diagonal, 
and then using these diagonal 1s to eliminate the nonzero entries in the
corresponding row and column.
Refer to reference \cite[\S 8.4.2]{BD-book} for a detailed discussion of the partial Smith form.

%%%%%%%%%%%%%%%%%%%%%%%%%%%%%%%%%%%%%%%%%%%%%%%%%%%%%%%%%%%%%%%%%%%%%%%%%%%%%%%%%%%%%%%%%%%%%%%%

\begin{table}
\begin{tabular}{lrrrrrrrr}
& rank & generators & nonzero & $\approx$\%nz & monic & degrees & Gbasis & radical
\\
Case 1
&\quad 1 & 340 & 102 & 27.4 & 93 & 1--4 & 9 & yes
\\
&\quad 2 & 21420 & 3263 & 15.2 & 2518 & 2--6& 45  & no
\\
&\quad 3 & 395080 & 52925 & 13.4 & 38305 & 3--8& 165  & no
\\
&\quad 4 & 2024785 & 380377 & 18.8 & 271515 & 4--10 & 548 & no
\\
Case 4
&\quad 1 & 216 & 38 & 17.6 & 30 & 1--3 & 6 & yes
\\
&\quad 2 & 8586 & 496 & 5.78 & 306 & 2--5 & 21 & no
\\
&\quad 3 & 99216 & 3053 & 3.08 & 1828 & 3--7 & 56 & no
\\
&\quad 4 & 316251 & 12791 & 4.04 & 7799 & 4--8 & 126 & no
\\
\end{tabular}
\smallskip
\caption{Gr\"obner bases for determinantal ideals in Cases 1 and 4}
\label{Grobnerinfo}
\end{table}

\subsubsection{First determinantal ideal}

Information about the computation of Gr\"obner bases for the determinantal ideals
is summarized in Table \ref{Grobnerinfo}.
For Case 1, the first determinantal ideal has $85 \cdot 4 = 340$ generators, 
only 102 of which are nonzero ($\approx 27.4\%$),
and after making these nonzero generators monic, we are left with 93 distinct generators
of degrees $1,\dots,4$.
The deglex Gr\"obner basis for this ideal consists of 9 polynomials:
\[
a_2, \quad 
a_3, \quad
a_5, \quad
a_6, \quad
a_7, \quad
a_8, \quad
a_9, \quad
a_4(a_4+1), \quad
a_{10}(a_{10}+1).
\]
From this it can be seen that this is a radical ideal.
The zero set consists of 4 distinct coefficient vectors
$a_1, \dots, a_{10}$ for which the matrix $M$ has rank 16:
\begin{equation}
\label{case1rank1solutions}
\left[
\begin{array}{rrrrrrrrrr}
 1 &\;   0 &\;   0 &\;    0 &\;   0 &\;   0 &\;   0 &\;   0 &\;   0 &\;   0  \\
 1 &\;   0 &\;   0 &\;    0 &\;   0 &\;   0 &\;   0 &\;   0 &\;   0 &\;  -1  \\
 1 &\;   0 &\;   0 &\;   -1 &\;   0 &\;   0 &\;   0 &\;   0 &\;   0 &\;   0  \\
 1 &\;   0 &\;   0 &\;   -1 &\;   0 &\;   0 &\;   0 &\;   0 &\;   0 &\;  -1  
\end{array}
\right]
\end{equation}
These coefficient vectors correspond to the following 4 operator identities:
\begin{equation}
\label{rank16identities1}
\begin{array}{l}
L^3(xy) = 0, \qquad
L^3(xy) = xL^3(y), \qquad
L^3(xy) = L^3(x)y,
\\[3pt]
L^3(xy) = L^3(x)y + xL^3(y).
\end{array}
\end{equation}
The last identity states that $L^3$ is a derivation.

\subsubsection{Second and third determinantal ideals}

In both cases, the radical of the ideal equals the first determinantal ideal,
and hence the zero sets are the same.
Hence there are no operator identities giving rank 17 or 18 of the matrix
of consequences.

\subsubsection{Fourth determinantal ideal}

The deglex Gr\"obner basis of the radical consists of these 19 polynomials:
\begin{align*}
&
a_6, \qquad 
a_7, \qquad 
a_9, \qquad 
a_2^2 - a_3, \qquad 
a_2 a_5, \qquad 
a_8 a_2, \qquad 
a_{10} a_2, \qquad 
-a_2 a_4 + a_3^2, 
\\
&
a_3 a_4 - a_2, \qquad
a_3 a_5, \qquad 
a_8 a_3, \qquad 
a_{10} a_3, \qquad 
-a_2 a_3 + a_4^2 - a_3 + a_4, \qquad 
a_4 a_5,
\\
&
a_8 a_4, \qquad 
a_5^2 - a_8, \qquad 
-a_{10} a_5 + a_8^2, \qquad 
a_{10} a_8 - a_5, \qquad 
a_{10}^2 - a_5 a_8 + a_{10} - a_8.
\end{align*}
The zero set of this ideal consists of the 4 points in display \eqref{case1rank1solutions}
together with the following 8 new points ($i = \sqrt{-1})$
for which the matrix $M$ has rank 19:
\begin{equation}
\label{case1rank4solutions}
\left[
\begin{array}{rrrrrrrrrr}
1 &\;   0 &\;   0 &\;   0 &\;   1 &\;  0 &\;  0 &\;   1 &\;  0 &\;   1 \\
1 &\;   0 &\;   0 &\;   0 &\;  -1 &\;  0 &\;  0 &\;   1 &\;  0 &\;  -1 \\ 
1 &\;   0 &\;   0 &\;   0 &\;   i &\;  0 &\;  0 &\;  -1 &\;  0 &\;  -i \\
1 &\;   0 &\;   0 &\;   0 &\;  -i &\;  0 &\;  0 &\;  -1 &\;  0 &\;   i \\
1 &\;   1 &\;   1 &\;   1 &\;   0 &\;  0 &\;  0 &\;   0 &\;  0 &\;   0 \\
1 &\;  -1 &\;   1 &\;  -1 &\;   0 &\;  0 &\;  0 &\;   0 &\;  0 &\;   0 \\
1 &\;   i &\;  -1 &\;  -i &\;   0 &\;  0 &\;  0 &\;   0 &\;  0 &\;   0 \\
1 &\;  -i &\;  -1 &\;   i &\;   0 &\;  0 &\;  0 &\;   0 &\;  0 &\;   0
\end{array}
\right]
\end{equation}
These coefficient vectors correspond to the following 8 operator identities:
\begin{equation}
\label{rank19identities}
\begin{array}{l}
   L^3(xy)  %   1
+  L^2(xL(y))  %   5 
+  L(xL^2(y))  %   8 
+  xL^3(y)  %  10 
= 0
\\[2pt]
   L^3(xy)  %   1
-  L^2(xL(y))  %   5 
+  L(xL^2(y))  %   8 
-  xL^3(y)  %  10 
= 0
\\[2pt]
    L^3(xy)  %   1
+i\,  L^2(xL(y))  %   5 
-   L(xL^2(y))  %   8 
-i\,  xL^3(y)  %  10 
= 0
\\[2pt]
    L^3(xy)  %   1
-i\,  L^2(xL(y))  %   5 
-   L(xL^2(y))  %   8 
+i\,  xL^3(y)  %  10 
= 0
\\[2pt]
   L^3(xy)  %   1
+  L^2(L(x)y)  %   2 
+  L(L^2(x))y  %   3 
+  L^3(x)y  %   4 
= 0
\\[2pt]
   L^3(xy)  %    1
-  L^2(L(x)y)  %    2 
+  L(L^2(x))y  %    3 
-  L^3(x)y  %    4 
= 0
\\[2pt]
    L^3(xy)  %    1
+i\,  L^2(L(x)y)  %    2 
-   L(L^2(x))y  %    3 
-i\,  L^3(x)y  %    4 
= 0
\\[2pt]
    L^3(xy)  %    1
-i\,  L^2(L(x)y)  %    2 
-   L(L^2(x))y  %    3 
+i\,  L^3(x)y  %    4 
= 0
\end{array}
\end{equation}

%%%%%%%%%%%%%%%%%%%%%%%%%%%%%%%%%%%%%%%%%%%%%%%%%%%%%%%%%%%%%%%%%%%%%%%%%%%%%%%%%%%%%%%%%%%%%%%%

\subsection{Case 2: $a_1 = 0$, $a_2 = 1$, and $a_3, \dots, a_{10}$ are free}

With the given substitutions, the matrix $M$ has full rank 20:
it does not have submaximal rank for any values of the parameters.
Hence there are no new operator identities in this case.

%%%%%%%%%%%%%%%%%%%%%%%%%%%%%%%%%%%%%%%%%%%%%%%%%%%%%%%%%%%%%%%%%%%%%%%%%%%%%%%%%%%%%%%%%%%%%%%%

\subsection{Case 3: $a_1, a_2 = 0$, $a_3 = 1$, and $a_4, \dots, a_{10}$ are free}

Similar to Case 2.

%%%%%%%%%%%%%%%%%%%%%%%%%%%%%%%%%%%%%%%%%%%%%%%%%%%%%%%%%%%%%%%%%%%%%%%%%%%%%%%%%%%%%%%%%%%%%%%%

\subsection{Case 4: $a_1, a_2, a_3 = 0$, $a_4 = 1$, and $a_5, \dots, a_{10}$ are free}

We make the indicated substitutions into the matrix $M$ of consequences.
We compute the partial Smith form of the resulting matrix and obtain
\[
M
\xrightarrow{\quad\text{PSF}\quad}
\left[
\begin{array}{cc}
I_{16} &\!\! O_{16,4} \\
O_{89,16} &\!\! C
\end{array}
\right]
\]
where the lower right block $C$ has size $89 \times 4$.
We delete the 35 zero rows in $C$ leaving a $54 \times 4$ matrix also called $C$.

\subsubsection{First determinantal ideal}

The deglex Gr\"obner basis consists the variables $a_5, \dots, a_{10}$
from which we see that this is a radical ideal.
The zero set is
\begin{equation}
\label{case1rank1}
\left[
\begin{array}{rrrrrrrrrr}
 0 &\;   0 &\;   0 &\;    1 &\;   0 &\;   0 &\;   0 &\;   0 &\;   0 &\;   0 
\end{array}
\right]
\end{equation}
which gives the following operator identity of rank 16:
\begin{equation}
\label{rank16identities2}
\begin{array}{l}
L^3(x)y = 0.
\end{array}
\end{equation}

\subsubsection{Second, third, and fourth determinantal ideals}

In all three cases, the radical ideal equals the first determinantal ideal.
Hence in Case 4, there are no operator identities which produce ranks 17, 18, or 19.

%%%%%%%%%%%%%%%%%%%%%%%%%%%%%%%%%%%%%%%%%%%%%%%%%%%%%%%%%%%%%%%%%%%%%%%%%%%%%%%%%%%%%%%%%%%%%%%%

\subsection{Cases 5 to 9} Similar to Case 2.

%%%%%%%%%%%%%%%%%%%%%%%%%%%%%%%%%%%%%%%%%%%%%%%%%%%%%%%%%%%%%%%%%%%%%%%%%%%%%%%%%%%%%%%%%%%%%%%%

\subsection{Case 10: $a_1, \dots, a_9 = 0$, $a_{10} = 1$}

In this case the matrix $M$ is a scalar matrix of rank 16.
We obtain the new operator identity
\begin{equation}
\label{rank16identities3}
\begin{array}{l}
x L^3(y) = 0.
\end{array}
\end{equation}

%%%%%%%%%%%%%%%%%%%%%%%%%%%%%%%%%%%%%%%%%%%%%%%%%%%%%%%%%%%%%%%%%%%%%%%%%%%%%%%%%%%%%%%%%%%%%%%%
%%%%%%%%%%%%%%%%%%%%%%%%%%%%%%%%%%%%%%%%%%%%%%%%%%%%%%%%%%%%%%%%%%%%%%%%%%%%%%%%%%%%%%%%%%%%%%%%
%%%%%%%%%%%%%%%%%%%%%%%%%%%%%%%%%%%%%%%%%%%%%%%%%%%%%%%%%%%%%%%%%%%%%%%%%%%%%%%%%%%%%%%%%%%%%%%%

\section{Main theorem}

The ranks of all the operator identities that we found were checked by
substitution of the coefficients into the matrix of consequences
and computing the rank.

\begin{theorem}
\label{theorem23}
For degree 2 and multiplicity 3,
the matrix of consequences has rank 16 if and only if the coefficients 
$a_1, \dots a_{10}$ correspond to the operator identities
in displays \eqref{rank16identities1}, \eqref{rank16identities2}, \eqref{rank16identities3}.
The matrix has rank 19 if and only if the coefficients correspond to the operator identities
in display \eqref{rank19identities}.
For all other values of the coefficients, the matrix has full rank 20.
\end{theorem}

%%%%%%%%%%%%%%%%%%%%%%%%%%%%%%%%%%%%%%%%%%%%%%%%%%%%%%%%%%%%%%%%%%%%%%%%%%%%%%%%%%%%%%%%%%%%%%%%
%%%%%%%%%%%%%%%%%%%%%%%%%%%%%%%%%%%%%%%%%%%%%%%%%%%%%%%%%%%%%%%%%%%%%%%%%%%%%%%%%%%%%%%%%%%%%%%%
%%%%%%%%%%%%%%%%%%%%%%%%%%%%%%%%%%%%%%%%%%%%%%%%%%%%%%%%%%%%%%%%%%%%%%%%%%%%%%%%%%%%%%%%%%%%%%%%

%%%%%%%%%%%%%%%%%%%%%%%%%%%%%%%%%%%%%%%%%%%%%%%%%%%%%%%%%%%%%%%%%%%%%%%%%%%%%%%%%%%%%%%%%%%%%%%%
%%%%%%%%%%%%%%%%%%%%%%%%%%%%%%%%%%%%%%%%%%%%%%%%%%%%%%%%%%%%%%%%%%%%%%%%%%%%%%%%%%%%%%%%%%%%%%%%
%%%%%%%%%%%%%%%%%%%%%%%%%%%%%%%%%%%%%%%%%%%%%%%%%%%%%%%%%%%%%%%%%%%%%%%%%%%%%%%%%%%%%%%%%%%%%%%%

\end{document}